\documentclass[12]{amsart}

\usepackage[all]{xy}
\usepackage{amssymb}
\usepackage{graphicx}

\newtheorem{thm}{Theorem}[section]
\newtheorem{prop}[thm]{Proposition}
\newtheorem{cor}[thm]{Corollary}
\newtheorem{lem}[thm]{Lemma}
\newtheorem{dfn}[thm]{Definition}
\newtheorem{example}[thm]{Example}

\newtheorem{remark}[thm]{Remark}

\def\e{\mathbf{c}}
\def\l{\mathbf{l}}
\def\r{\mathbf{r}}

\def\t{\mathbf{t}}
\def\d{\mathbf{d}}

\newcommand{\F}{{\check{\pi}}_0}
\def\bF{{\check{\bar{\pi}}_0}}
\def\C{\check{c}}
\def\bC{\check{\bar{c}}}

\newcommand{\CC}{\mathbb{C}}
\def\R{\mathbb{R}}
\def \N {\mathbb{N}}
\def \id {\mathrm{id}}
\def\tr{\mathrm tr}
\def \den{{\rm den}}
\def \inv{{\rm inv}}
\def \Fr {\mathrm{Fr}}
\newcommand{\Z}{\mathbb{Z}}

\begin{document}

\title[Limit and end functors of dynam. systems via
exterior spaces]{Limit and end functors of dynamical systems via
exterior spaces}

\thanks{Partially supported by Ministerio de
Educaci\'on y Ciencia grant MTM2009-12081, Universty of La Rioja project API12/10  and FEDER}



\author{J.M. Garc\'{\i}a Calcines, L. Hern\'andez Paricio and M. Teresa Rivas
Rodr\'{\i}guez}

\address{J.M. Garc\'{\i}a Calcines \newline \indent Departamento de Matem\'{a}tica
Fundamental \newline \indent Universidad de La Laguna \newline
\indent 38271 LA LAGUNA. } \email{jmgarcal@ull.es}

\address{L. Hern\'andez Paricio; \hspace{6pt}
M. Teresa Rivas Rodr\'{\i}guez \newline \indent Departamento de
Matem\'{a}ticas y Computaci\'{o}n \newline \indent  Universidad de
La Rioja \newline \indent 26004 LOGRO\~{N}O.}
\email{luis-javier.hernandez@unirioja.es}
\email{maria-teresa.rivas@unirioja.es}



\begin{abstract}
In this paper we analyze some applications of the category of
exterior spaces to the study of dynamical systems (flows). We
study the notion of an absorbing open subset of a dynamical
system; i.e., an open subset that contains the ``future part'' of
all the trajectories. The family of all absorbing open subsets is
a quasi-filter  which gives  the structure of an exterior space to
the flow. The  limit space and end space of  an exterior space are
used to construct the limit spaces and end spaces of a dynamical
system. On the one hand, for a dynamical system two limits spaces
$L^{\r}(X)$ and $\bar L^{\r}(X)$ are constructed  and their
relations with the subflows of periodic, Poisson stable points and
$\Omega^{\r}$-limits of $X$ are analyzed. On the other hand,
different end spaces are also associated to a dynamical system
having the property that any positive semi-trajectory has an end
point in these end spaces. This type of  construction permits us
to consider the subflow containing all  trajectories finishing at
an end point $a$. When $a$ runs over the set of all end points, we
have an induced decomposition of a dynamical system as a disjoint
union of stable (at infinity) subflows.
\end{abstract}

\maketitle

\section{Introduction}

Many natural phenomena can be modelized by means of a family of
differential equations and each one can be put (maybe after some
manipulations) in the form $\dot \phi =f(\phi)$, where $\phi$ are
local  coordinates in an open neighborhood of the  point $p \in
M$, where $M$ is an $m$-dimensional manifold, $\dot \phi $ are the
coordinates of tangent vectors in the open neighborhood of the
point $p\in M$ and $f$ is a real valued function whose domain is
an open of $\R^m$.

Under the assumption of $f$ being locally lipschitzian, an initial
condition $\phi^p (0) = p$ uniquely determines a maximal solution
$\phi^p (t)$. However, the domain of $\phi^p (t)$ does not need to
be the whole real line $\R$, but only an open interval, $(a^p,b^p)
\subset \R \,, \, a^p<0<b^p$,  which depends on the initial
condition. All the solutions give  a  local flow $\phi \colon W
\to M$, $\phi(t,p)=\phi^p (t)$, where  $W$  is an open subset of
$\R \times M $ containing $ \{0\}\times M$ and if we denote
$\phi_s(p)=\phi(s,p)$, $(s,p)\in W$, $\phi$ satisfies $\phi_0
=\id_M$, $\phi_t \phi_s = \phi_{t+s}$, wherever it makes sense.
The space $M$ is called the  phase space and $\phi$ is also called
the phase map. The trajectory of a point $p\in M$ is the subset
$\gamma(p)=\{\phi(t,p)| t\in (a^p,b^p) \}$. It is easy to check
that $M$ is a disjoint union of trajectories. We note that when a
trajectory has more than one point, a natural orientation is
induced by the canonical orientation of the real numbers $\R$.
Then, we can consider $M$ as a disjoint union of critical
trajectories and oriented trajectories to obtain a phase portrait
of the dynamical system $\phi$. It is well known that (under some
assumptions, see \cite{Bhatia}) if  $\varphi$ is a local flow on
$M$, then there exists a global ($W= \R\times M$) flow $\phi$ in
$M$ such that the oriented trajectories of $\varphi$ and $\phi$
coincide. Consequently their phase portraits are the same. As a
consequence of this type of result we have considered convenient
to reduce our study to the case of global flows. This allows us a
functorial approach to the study of dynamical systems and their
properties.

On the other side, for every topological space $X$, a continuous
action $\varphi \colon \R \times X \to X$ induces a group
homomorphism $\bar \varphi \colon \R \to \mbox{Aut}(X)$, where $
\mbox{Aut}(X)$ is the group of homeomorphisms of $X$ provided with
the compact-open topology. This allows us to study a flow as a
particular case of a transformation group.

The pioneering work of H. Poincar{\'e} \cite{Poincare,Poincare2}
in the late XIX century studied the topological properties of
solutions of autonomous ordinary differential equations. We can
also mention the work of A. M. Lyapunov \cite{Liapunov} who
developed his theory of stability of a motion (solution) of a
system of $n$ first order ordinary differential equations. While
much of Poincar{\'e}'s work analyzed the global properties of the
system, Lyapunov's work looks at the local stability of a
dynamical system. The theory of dynamical systems reached a great
development  with the work of G.D. Birkhoff  \cite{Birkhoff}, who
may be considered as the founder of this theory. He established
two main lines in the study of dynamical systems: the topological
theory and the ergodic theory.

Previously, the
authors have developed some results  on proper homotopy theory and
exterior homotopy theory to classify non-compact spaces and to
study the shape of a compact metric space. 
In this paper, we describe some basic ideas that permit a new
treatment of the study of dynamical system. The main
objective proposed by the authors is to describe functorially a
way of using exterior spaces to study dynamical systems. In our
approach, the main key to establish a connection between exterior
spaces and dynamical systems is the notion of absorbing open
region (or $\r$-exterior subset). Given a flow on space $X$,  an
open set $E$ is said to be  $\r$-exterior if for every $x \in X$
there is $r_0 \in \R$ such that $r\cdot x \in E$ for $r\geq r_0$.
The space $X$ together with the family of $\r$-exterior open
subsets is an exterior space which is denoted by $X^{\r}$.

On the one hand, using exterior spaces, we construct two limit
subflows $L^{\r}(X)$ and $\bar L^{\r}(X)$ associated with a flow
$X$. One of the important results of our study (see Corollary
\ref{reticulo}) is the following:

{\it If $X$ is a locally compact $T_3$ space, then $ L^{\r}(X) =
P(X)$ and $\bar L^{\r}(X) = \overline{\Omega^{\r}(X)}$;
furthermore we have that}
$$L^{\r}(X)=P(X)\subset P^{\r}(X)\subset \Omega^{\r}(X) \subset
\overline{\Omega^{\r}(X)}=\bar L^{\r}(X),$$ where $P(X)$ is the
subflow of periodic points, $P^{\r}(X)$ is the subflow of right
(positive) Poisson stable points and $\Omega^{\r}(X)=
\bigcup_{x\in X} \omega^{\r}(x)$ ($\omega^{\r}(x)$ is the right
limit set of the point $x\in X$). Under the condition of being
$\r$-regular at infinity (see Proposition \ref{eleigualbarele}),
we have that $L^{\r}(X)=\bar L^{\r}(X)$ and in some cases we can
also ensure that $\bar L^{\r}(X)$ is compact. The constructions
$L^{\r}, \bar L^{\r}$ induce  classifications (for flows) of the
following type: Two flows $X,Y$ are said to be $\bar L$-equivalent
if there is a flow morphism $f\colon X \to Y$ such that $\bar L
(f)$ is a homotopy equivalence. This will determine equivalence
classes  of flows having, for instance,  the same number of
critical points or the `same type' of periodic trajectories.

On the other hand, using again exterior spaces, we construct
`slightly different' end spaces $\F^{\r}(X)$,  $\C^{\r}(X)$,
$\bF^{\r}(X)$, $\bC^{\r}(X)$, which agree when $X$ is locally
path-connected and ${\r}$-regular at infinity. In this case, the
end spaces have a pro-discrete topology (with additional
conditions, a pro-finite topology). The importance of this end
space is that each right semi-trajectory of the flow has an end
point in this space. This fact allows one to give a set map $\bar
\omega_{\r} \colon X \to \bC^{\r}(X)$ and the corresponding $\bar
\omega_{\r}$-decomposition of the flow
$$ X=\bigsqcup_{a\in \bC^{\r}(X) }\bar X_{(\r, a)}.$$
In general, there will be end points which are not reached by  right
semi-trajectories and  end points of semi-trajectories
that are not reached by  right semi-trajectories contained in the
limit subflow. In this paper, we give sufficient conditions to
ensure that an end point can be reached  by right
semi-trajectories of $\bar L(X)$ (see Proposition
\ref{erepresentableuno}). These end spaces will be used (not in
this work) to construct completions (compactification under some
topological and dynamical conditions) of a flow. These completions
are related to Freundenthal compactifications
\cite{EHR93,Freudenthal, Kererjarto} and will permit us to apply some
nice properties of compact flows to a more general class of
topological flows.

It is important to observe that applying the results of this paper
to the reversed flow we will obtain all the corresponding dual
concepts, constructions and properties.

\section{Preliminaries on exterior spaces and dynamical systems}

\subsection{The categories of proper and exterior spaces}\label{exteriorspaces}

A continuous map $f:X\rightarrow Y$ is said to be proper if for
every closed compact subset $K$ of $Y$, $f^{-1}(K)$ is a compact
subset of $X$. The category of topological  spaces and the
subcategory of spaces and proper maps will be denoted by {\bf Top}
and {\bf P}, respectively. This last category and its
corresponding proper homotopy category are very useful for the
study of non-compact spaces. Nevertheless, one has the problem
that ${\bf P}$ does not have enough limits and colimits and then
we can not develop the usual homotopy constructions such as loops,
homotopy limits and colimits, et cetera. An answer to this problem
is given by  the notion of exterior space. The new  category of
exterior spaces and maps is complete and cocomplete and contains
as a full subcategory the category of spaces and proper maps, see
\cite{GGH98, GGH04}. For more properties and applications of
exterior homotopy categories we refer the reader to \cite{GGH01,
EHR05, DHR09, GHR09, H95} and for  a survey of proper homotopy to
\cite{P95}.

\begin{dfn}
Let $(X,\t )$ be a topological space, where $X$ is the subjacent
set and $\t$ its topology. An \emph{externology} on $(X,\t )$ is a
non-empty collection $\varepsilon $ (also denoted by $\varepsilon
(X)$) of open subsets which is closed under finite intersections
and such that if $E\in \varepsilon $, $U\in \t$ and $E\subset U$
then $U\in \varepsilon .$ The members of $\varepsilon $ are called
exterior open subsets. An \emph{exterior space} $(X,\varepsilon,
\t )$ consists of a space $(X,\t )$ together with an externology
$\varepsilon $. A map $f:(X,\varepsilon, \t)\rightarrow
(X',\varepsilon', \t')$ is said to be an \emph{exterior map} if it
is continuous and $f^{-1}(E)\in \varepsilon $, for all $E\in
\varepsilon '.$
\end{dfn}

The category of exterior spaces and maps will be denoted by {\bf
E}. Given an space $(X, \t_{X})$, we can always consider the
trivial exterior space taking $\varepsilon = \{ X \}$ or the total
exterior space if one takes $ \varepsilon = \t_{X}$. An important
example of externology on a given topological space $X$ is the one
constituted by the complements of all closed-compact subsets of
$X$, that will be called the cocompact externology and usually
written as $\varepsilon ^{\e}(X).$ The category of spaces and
proper maps can be considered as a full subcategory  of the
category of exterior spaces via the full embedding
$(\cdot)^{\e}:{\bf P}\hookrightarrow {\bf E}.$ The functor
$(\cdot)^{\e}$ carries a space $X$ to the exterior space $X^{\e}$
which is provided with the topology of $X$ and the externology
$\varepsilon ^{\e}(X)$. A map $f\colon X \to Y$ is carried to the
exterior map $f^{\e}\colon X^{\e} \to Y^{\e}$ given by $f^{\e} =
f$. It is easy to check that a continuous map $f:X\rightarrow Y $
is proper if and only if $f=f^{\e}\colon X^{\e} \to Y^{\e}$  is
exterior.

An important role in this paper will be played by the following
construction $(\cdot)\bar \times (\cdot) $: Let $(X, \varepsilon
(X),\t_X)$ be an exterior space, $(Y,\t_Y)$ a topological space
and for $y\in Y$ we denote by $(\t_{Y})_y$ the family of open
neighborhoods of $Y$ at $y$. We consider on $X\times Y$ the
product topology $\t_{X \times Y}$ and the externology
$\varepsilon ({X\bar \times Y})$ given by those $E\in \t_{X\times
Y}$ such that for each $y\in Y$ there exists
 $U_y\in (\t_{Y})_y$ and $T^{y}\in \varepsilon ({X})$ such that
$T^{y}\times U_{y}\subset E.$ This exterior space will be denoted
by $X\bar{\times }Y$ in order to avoid a possible confusion with
the product externology. This construction gives a functor
$$(\cdot)\bar \times (\cdot) \colon {\bf E}\times {\bf Top}\rightarrow
{\bf E}.$$ When $Y$ is a compact space,  we have that $E$ is an
exterior open subset if and only if it is an open subset and there
exists $G\in \varepsilon ({X})$ such that $G\times Y\subset E$.
Furthermore, if $Y$ is a compact space and $\varepsilon (X)
=\varepsilon ^{\e}(X),$ then $\varepsilon (X\bar{\times } Y)$
coincides with $\varepsilon ^c(X\times Y)$ the externology of the
complements of closed-compact subsets of $X\times Y.$ We also note
that if $Y$ is a discrete space, then $E$ is an exterior open
subset if and only if it is open  and for each $y\in Y$ there is
$T^y\in \varepsilon ({X})$ such that $T^y\times \{y\}\subset E$.

This bar construction provides a natural way to define {\em
exterior homotopy} in ${\bf E}$. Indeed, if $I$ denotes the closed
unit interval, given exterior maps $f,g:X\rightarrow Y,$ it is
said that $f$ \emph{is exterior homotopic to} $g$ if there exists
an exterior map $H:X\bar{\times }I\rightarrow Y$ (called exterior
homotopy) such that $H(x,0)=f(x)$ and $H(x,1)=g(x),$ for all $x\in
X.$ The corresponding homotopy category of exterior spaces will be
denoted by $\pi {\bf E}.$ Similarly, the usual homotopy category
of topological spaces will be denoted by $\pi {\bf Top}.$

\subsection{Dynamical Systems and $\Omega$-Limits}

Next we recall some elementary concepts about dynamical systems.

\begin{dfn}
A {\sl dynamical system} (or a {\sl flow}) on a topological  space
$X$ is a continuous map $\varphi \colon \R {\times} X\to  X$,
$\varphi(t,x)=t \cdot x$,  such that
\begin{enumerate}
\item[(i)] $\varphi(0,p)=p,$ $\forall p \in X;$
\item[(ii)] $\varphi(t,\varphi(s, p))=\varphi( t+s,p),$ $\forall p
\in X,$ $\forall t,s \in \mathbb{R}.$
\end{enumerate}
A flow on $X$ will be denoted by $(X, \varphi)$ and when no
confusion is possible, we use $X$ for short.
\end{dfn}

For a subset $A \subset X$, we denote inv$(A)=\{p \in A| \R \cdot
p \subset A\}$.

\begin{dfn}
A subset $S$ of a flow $X$ is said to be {\sl invariant} if
inv$(S)=S.$
\end{dfn}

Given a flow $\varphi \colon \R {\times} X\to  X$ one has a
subgroup $\{\varphi_t \colon X \to X| t\in \R\}$ of
homeomorphisms, $\varphi_t(x) = \varphi(t,x)$, and a family of
motions  $\{\varphi^p\colon \R \to X| p \in X\}$, $\varphi^p(t) =
\varphi(t,p).$

\begin{dfn} Given two flows  $\varphi \colon \R {\times} X\to  X$,
$\psi \colon \R {\times} Y\to  Y$, a  \emph{flow morphism}
$f\colon (X, \varphi) \to (Y, \psi)$ is a continuous map $f\colon
X \to Y$ such that $f(r\cdot p)=r\cdot f(p)$ for every $r\in \R$
and for every $p\in X$.
\end{dfn}

We note that if $S \subset X$ is invariant, $S$ has a flow
structure and the inclusion is a flow morphism.

We denote by $\bf{ F}$ the category of flows and flows morphisms.

\medskip

We recall some basic fundamental examples: (1) $X=\R$ with the
action $\varphi \colon \R \times X \to X$,  $\varphi(r, s)= r+s$.
(2) $X=S^1=\{z\in {\mathbb{C}}| |z|=1\}$ with  $\varphi \colon \R
\times X \to X$, $\varphi(r, z)= e^{2\pi i r}z$. (3) $X=\{0\}$
with the trivial action $\varphi \colon \R \times X \to X$ given
by $\varphi(r, 0)= 0.$ In all these cases, the flows only have one
trajectory.

\begin{dfn}\label{rightlimit} For a flow $X$,
the {\sl $\omega^{\r}$-limit set} (or right-limit set, or positive
limit set) of a point $p \in X$ is given as follows:
$$\omega^{\r} (p) =\{q \in X | \exists \hspace{3pt}\mbox{a net}\hspace{3pt}t_{\delta } \to +\infty \mbox{ such that }
t_{\delta } \cdot p \to q\}.$$
\end{dfn}

 If $\overline{A}$ denotes the closure of a subset $A$ of a topological space,
we  note that the subset $\omega^{\r}(p)$ admits the alternative
definition
$$\omega^{\r} (p) = \bigcap_{t \geq 0}\overline{ [t, +\infty) \cdot p}$$
which has the advantage of showing that $\omega^{\r}(p)$ is
closed.

\begin{dfn}
The {\sl $\Omega^{\r}$-limit set} of a flow $X$  is the following
invariant subset:
$$\Omega^{\r} (X) =\bigcup_{p\in X} \omega^{\r}(p)$$
\end{dfn}

Now we introduce the basic notions of critical, periodic and
$\r$-Poisson stable points.

\begin{dfn} Let $X$ be a flow.
A point $x\in X$ is said to be  a {\sl critical point} (or a
\emph{rest point}, or an \emph{equilibrium point}) if for every
$r\in \R,$ $r\cdot x =x.$ We denote by $C(X)$ the invariant subset
of critical points of $X.$
\end{dfn}

\begin{dfn} Let $X$ be a flow. A point $x\in X$ is said to be
\emph{periodic} if there is $r\in \R$, $r\not = 0$ such that
$r\cdot x= x.$ We denote by $P(X)$ the invariant subset of
periodic points of $X.$
\end{dfn}

It is clear that a critical point is a periodic point. Then
$$C(X)\subset P(X).$$

If $x \in X$ is a periodic point but not critical, then the real
$r\not =0$ such that $r\cdot x=x$ and $r$ is called a
\emph{period} of $x.$ The smallest positive period $r_0$ of $x$ is
called the {\it fundamental period} of $x.$ Further if $r \in \R$
is such that $r\cdot x=x,$ then there is $z\in \Z$ such that
$r=zr_0.$

\begin{dfn}\label{poisson} Let $(X, \varphi)$ be a flow. A point $x \in X$
is said to be $\r$-\emph{Poisson stable} if there is a net
$t_{\delta }\to +\infty$ such that $t_{\delta } \cdot x \to x;$
that is, $x\in \omega^{\r}(x).$ We will denote by $P^{\r}(X)$ the
invariant subset of $\r$-Poisson stable points of $X.$
\end{dfn}

The reader can easily check that $$P(X)\subset P^{\r}(X) \subset
\Omega^{\r}(X).$$

 The notions above can be dualized to obtain the notion of the
$\omega^{\l}$-limit ($\l$ for `left') set of a point $p$, the
$\Omega^{\l}$-{ limit} of  $X$, $\l$-Poison stable points, et
cetera.

\begin{remark}
Observe that when $X$ satisfies the first axiom of countability
(for instance, when $X$ is metrizable) we can consider sequences
instead of nets in definitions \ref{rightlimit} and \ref{poisson}.
\end{remark}

\section{End Spaces and Limit Spaces of an exterior space}\label{tres}

In this section we will deal with special limit constructions.
Observe that if  $X=(X, \varepsilon (X))$ is an exterior space,
then its externology $\varepsilon(X)$ can be seen as an inverse
system of spaces. Composing with different endofunctors
$\mathbf{Top}\rightarrow \mathbf{Top}$ and taking the
corresponding topological inverse limit we will be able to obtain
functors $\mathbf{E}\rightarrow \mathbf{Top}$.

\subsection{The functors $L, \F, \C \colon {\bf E} \to {\bf Top}$}

Given an exterior space $X=(X, \varepsilon (X))$, its externology
$\varepsilon(X)$ is an inverse system of spaces. Then we define
the limit space of $(X, \varepsilon (X))$ as the topological space
$$L(X)=\lim  \varepsilon (X).$$  Note that for each $E' \in
\varepsilon (X)$  the canonical map $\lim \varepsilon (X) \to E' $
is continuous and factorizes as $\lim \varepsilon (X) \to \cap
_{E\in \varepsilon (X)}E\to  E'.$ Therefore the canonical map
$\lim  \varepsilon (X) \to \cap _{E\in \varepsilon (X)}E$ is
continuous. On the other side, by the universal property of the
limit, the family of maps $\cap _{E\in \varepsilon
(X)}E\to E',$ $E' \in \varepsilon (X)$ induces a continuous map
$\cap _{E\in \varepsilon (X)}E\to \lim \varepsilon (X)$. This
implies that the canonical map $\lim \varepsilon(X) \to \cap_{E\in
\varepsilon(X)}E$ defines a natural homeomorphism.

We recall that for a topological space $Y$, $\pi_0(Y)$ denotes the
set of path-components of $Y$ and we have a  canonical map $Y\to
\pi_0(Y)$ which induces a quotient topology on $\pi_0(Y).$
Similarly, if $c(Y)$ denotes the set of connected components of a
space $Y$, we have a similar quotient map $Y\to c(Y).$ We remark
that if $Y$ is locally path-connected (respectively, locally
connected), then $\pi_0(Y)$ ($c(Y)$) is a discrete space.

It is also interesting to note that  for any topological space
$Y$, there exists  a canonical commutative diagram of natural
maps:
$$
\xymatrix{&Y\ar[ld]\ar[rd]& \\
\pi_0 (Y)\ar[rr]&&c(Y) }
$$

\begin{dfn} Given an exterior space $X=(X, \varepsilon(X))$
the \emph{limit space} of $X$ is the topological subspace
$$L(X)=\lim \varepsilon(X)= \cap_{E\in \varepsilon(X)}E.$$

The \emph{end space} of $X$ is the inverse limit
$$\F(X)=\lim \pi_0 \varepsilon(X)=\lim_{E\in \varepsilon(X)}\pi_0(E)$$
provided with the inverse limit topology of the spaces $\pi_0(E).$

The  $c$-end space of $X$ is the inverse limit
$$\C(X)=\lim c\, \varepsilon(X)=\lim_{E\in \varepsilon(X)}c(E)$$
provided with the inverse limit topology of the spaces $c(E).$ The
elements of $\F(X)$ or $\C(X)$ will be called \emph{end points} of
$X.$
\end{dfn}

An end point $a\in \F(X)$ is represented by the filter base
$$\{U_a^E | U_a^E \mbox{ is a path-component of }E,E \in \varepsilon (X)\}.$$

We note that a locally path-connected exterior space $(X,
\varepsilon(X))$ induces the following family of exterior spaces
$$\{ (X, \varepsilon(X,a))| a\in \F(X)\}$$
where $\varepsilon(X,a)$ is the externology generated by the
filter base
$$\{U_a^E | U_a^E \mbox{ is a path-component of }E, E \in \varepsilon (X)\}.$$
The end points of $\C (X)$ have similar properties.

It is interesting to observe that if $X$ is an exterior space and
$X$ is locally path-connected (locally connected), then $\F(X)$
($\C(X)$) is a prodiscrete space. On the other hand, given any
exterior space $(X, \varepsilon(X))$, we have a canonical
commutative diagram of natural maps
$$
\xymatrix{& L(X) \ar[ld]_{e_0} \ar[rd]^{e} & \\
\F (X) \ar[rr]_{\theta} & & \C (X) }
$$

\begin{dfn} Given an exterior space $X=(X, \varepsilon(X))$,
an end point $a\in \F(X)$ ($a\in \C(X)$) is said to be
$e_0$-\emph{representable} ($e$-\emph{representable}) if there is
$p\in L(X)$ such that  $e_0(p)=a$ ($e(p)=a$). Notice that the maps
$e_0\colon L(X)\to \F(X),$ $e\colon L(X) \to \C(X)$ induce an $
e_0$-decomposition and an $e$-decomposition
$$L(X)= \bigsqcup_{a \in  \F(X)} L_a^0(X),
\quad L(X)=\bigsqcup_{a \in  \C(X)} L_a(X)$$ where $L_a^0(X)=
e_0^{-1}(a)$ and $L_a(X)=e^{-1}(a)$. These special subsets will be
respectively called the $e_0$-component of the end $a\in \F(X)$
and the $e$-component of the end $a\in \C(X)$ in the limit $
L(X).$
\end{dfn}

We denote by $e_0L(X)$ and $eL(X)$ the corresponding subsets of
representable end points. It is clear that
$$ L(X)= \bigsqcup_{a \in e_0L(X)} L_a^0(X), \quad L(X)=
\bigsqcup_{a \in  eL(X)} L_a(X)$$ and for $b\in eL(X)$ one has
that$$ L_b(X)= \bigsqcup_{a \in ( \theta^{-1}(b) \cap e_0L(X))}
L_a^0(X).$$

\begin{example} \quad Let $M\colon \R \to (0,1)$ be an increasing
continuous map such that $\lim_{t\to -\infty}M(t)=0$ and
$\lim_{t\to +\infty}M(t)=1$ and take $A=\{e^{2\pi i t}|t\in \R\}$,
$B=\{M(t)e^{2\pi i t}|t\in \R\}$. Consider $X=A\cup B\subset \CC$
provided with the relative topology (observe that $X$ is not
locally connected). On the topological space $X$ the  flow
$\varphi \colon \R \times X \to X$ is given by $\varphi (r,e^{2\pi
i t})=e^{2\pi i (r+t)}$, $\varphi (r,M(t) e^{2\pi it})=M(r+t)
e^{2\pi i (r+t)}$. It is clear that this flow has two trajectories
$A,B$. If for each natural number $n$ we denote $B_n=\{M(t)e^{2\pi
i t}|t\geq n\}$, then a base of an externology on $X$  is given by
$\{E_n=A\cup B_n|n\in \N\}.$ Since $A,B_n$ are path-connected and
$\overline{B_n}=E_n$ is connected, it follows that
$\pi_0(E_n)=\{A,B_n\}$ and $c(E_n)=\{E_n\}$. Therefore
 $$\check \pi_0(X)=\{*_A, *_B\}, \quad \check c(X)=\{*\}$$
 For this example we have  $L(X)=A$, the $e_0$-decomposition
 $$L_{*_A}^0=A, \quad L_{*_B}^0=\emptyset$$
 and the $e$-decomposition $L_{*}=A$.
 This means  that $*_B $ is not $e_0$-representable.
\end{example}

 It is not difficult to check that the functor
$L$ preserves homotopies and the functors $\F,\C$ are invariant by
exterior homotopy.

\begin{lem}\label{lemahomotopia} Suppose that $X$ and $Y$ are exterior spaces and
$f,g \colon X\to Y$ exterior maps.
\begin{itemize}
\item[(i)] If $H \colon X\bar \times I \to Y$ is an exterior homotopy from $f$ to $g$,
then $L(H)=H|_{L(X) \times I} \colon  L(X\bar \times I)=L(X)\times
I  \to  L(Y)$ is a homotopy from $L(f)$ to $L(g).$
\item[(ii)] If $f$ is exterior homotopic to $g$, then $\F(f)=\F(g)$ and $\C(f)=\C(g).$
\end{itemize}
\end{lem}

As a consequence of this lemma one has:

\begin{prop}\label{homotopy} The functors $L, \F, \C\colon {\bf E}
\to {\bf Top}$ induce functors  $$L\colon \pi {\bf E} \to \pi{\bf
Top},\hspace{15pt}\F, \C \colon \pi{\bf E} \to {\bf Top}.$$
\end{prop}

It is interesting to observe that the functor $L\colon {\bf E} \to
{\bf Top}$ admits in a natural way a  left adjoint: Given a
topological space $X$, recall that we can consider on $X$ the
trivial externology $\varepsilon^{\tr}(X)=\{X\}.$ This
construction gives the exterior space $
X_{\tr}=(X,\varepsilon^{\tr}(X))$ and induces the canonical
functor $(\cdot)_{\tr} \colon {\bf Top} \to {\bf E}$, $X\mapsto
X_{\tr}$. The reader can straightforwardly check the following
result:

\begin{prop}\label{adjoint1} The functor $(\cdot)_{\tr} \colon {\bf Top}
\to {\bf E}$ is left adjoint to the functor $L\colon  {\bf E} \to
{\bf Top}.$  Moreover, this pair of adjoint functors induces an
adjunction on the homotopy categories: $(\cdot)_{\tr} \colon \pi
{\bf Top} \to \pi {\bf E}$,  $L\colon \pi {\bf E} \to \pi{\bf
Top}$.
\end{prop}

\subsection{The functors $\bar L, \bF, \bC \colon {\bf E} \to {\bf Top}$}

The  externology of an exterior space $X=(X, \varepsilon(X))$ and
the closure operator of the subjacent  topological space  induce
the following inverse system $\bar \varepsilon (X)=\{\overline{E}|
E \in  \varepsilon(X)\}.$ Using this new inverse system, we can
rewrite  notions and analogous results of subsection above as
follows:

\begin{dfn} Given an exterior space $X=(X, \varepsilon(X))$
the bar-limit space of $X$ is the topological subspace
$$\bar L(X)=\lim  \bar \varepsilon(X)= \cap_{E\in \varepsilon(X)}\overline{E}.$$

The  bar-end space of $X$ is the inverse limit
$$\bF(X)=\lim \pi_0
\bar  \varepsilon(X)=\lim_{E\in
\varepsilon(X)}\pi_0(\overline{E})$$ provided with the inverse
limit topology of the  spaces $\pi_0(\overline{E}).$

The  $c$-bar-end space of $X$ is the inverse limit
$$\bC(X)=\lim c\, \bar \varepsilon(X)=\lim_{E\in \varepsilon(X)}c(\overline{E})$$
provided with the inverse limit topology of the  spaces
$c(\overline{E}).$
\end{dfn}

Given any exterior space $X=(X, \varepsilon(X))$, we have a
canonical diagram of natural maps
$$
\xymatrix{&\bar L(X)\ar[ld]_{\bar e_0}\ar[rd]^{\bar e}& \\
\bF(X)\ar[rr]_{\bar \theta}&&\bC(X) }
$$
\noindent and there canonical natural maps $L(X) \subset
\overline{L(X)} \subset \bar L (X)$, $\F(X) \to  \bF(X)$, $\C(X)
\to  \bC(X)$ such that the following diagram is commutative:
$$
\xymatrix{
& L(X)\ar[ld] \ar[rd] \ar'[d][dd] & \\
\F(X) \ar[rr] \ar[dd] && \C(X)  \ar[dd] \\
& \bar L(X) \ar[ld] \ar[rd] & \\
 \bF(X) \ar[rr] && \bC(X) }
$$

\begin{dfn} Given an exterior space $X=(X, \varepsilon(X)),$
an end point $a\in \bF(X)$ ($a\in \bC(X)$) is said to be $\bar
e_0$-representable ($\bar e$-representable) if there is  $p\in
\bar L(X)$ such that  $ \bar e_0(p)=a$($\bar e(p)=a$). The maps $
\bar e_0 \colon  \bar L(X)\to \F(X), \bar e\colon \bar L(X) \to
\C(X)$ induce an $ \bar e_0$-decomposition and an $\bar
e$-decomposition
$$ \bar L(X)= \bigsqcup_{a \in  \bF(X)} \bar L_a^0(X), \quad \bar L(X)=
\bigsqcup_{a \in  \bC(X)} \bar L_a(X)$$ where $ \bar L_a^0(X)=
\bar e_0^{-1}(a)$ ($ \bar L_a(X)= \bar e^{-1}(a)$) will be called
the $\bar e_0$-component ($\bar e$-component) of the end $a\in
\bF(X)$ ($a\in \bC(X)$)  in the limit $ \bar L(X)$.
 \end{dfn}

 We denote by $\bar e_0\bar L(X)$ ($\bar e \bar L(X)$) the corresponding
 subset of representable end points. It is clear that we have a commutative diagram
$$\xymatrix{
& L(X)\ar[ld] \ar[rd] \ar'[d][dd] & \\
e_0L(X) \ar[rr] \ar[dd] && eL(X)  \ar[dd] \\
& \bar L(X) \ar[ld] \ar[rd] & \\
 \bar e_0\bar L(X) \ar[rr] && \bar e\bar L(X) }
$$
We also have the following similar results:

\begin{lem}\label{lemabarhomotopia} Suppose that $X, Y$ are exterior spaces and
$f, g \colon X\to Y$ exterior maps.
\begin{itemize}
\item[(i)] If $H \colon X\bar \times I \to Y$ is an exterior homotopy from $f$ to $g$,
then $\bar L(H)=H|_{\bar L(X) \times I}\colon  \bar L(X\bar \times
I)=\bar L(X)\times I \to \bar L(Y)$ is a homotopy from $\bar L(f)$
to $\bar L(g).$
\item[(ii)] If $f$ is exterior homotopic to $g$, then $\bF(f)=\bF(g)$ and
$\bC(f)=\bC(g).$
\end{itemize}
\end{lem}

As a consequence of this lemma we have:

\begin{prop}\label{homotopy}
The functors $\bar L, \bF, \bC\colon {\bf E} \to {\bf Top}$ induce
functors $$\bar L\colon \pi {\bf E} \to \pi{\bf
Top},\hspace{10pt}\bF, \bC \colon \pi{\bf E} \to {\bf Top}.$$
\end{prop}

In Proposition \ref{adjoint1} a left adjoint has been constructed
for the functor $L$. Nevertheless in the case of functor $\bar L$
we have the following alternative result:

\begin{prop} Suppose that $X, Y$ are exterior spaces and $X$
satisfies that for every $E \in \varepsilon (X)$,
$\overline{E}=X$. Then we have the following canonical injective
map $$Hom_{\bf E}(X,Y)\to Hom_{\bf Top}(X_{\t},\bar L(Y)),$$ where
$X_{\t}$ denote the subjacent topological space of $X.$ Therefore,
if ${\bf E}_{\den}$ denote the full subcategory of  exterior
spaces $X$ satisfying that for every $E \in \varepsilon (X)$,
$\overline{E}=X$, then $\bar L \colon{\bf E}_{\den} \to {\bf Top}$
is a faithful functor.
\end{prop}

\subsection{Topological and  exterior properties and canonical maps}

Let $X=(X, \varepsilon(X))$ be an exterior space and consider
$\bar\varepsilon(X)=\{\overline{E}| E\in \varepsilon(X)\}.$

\begin{dfn} An exterior space $X=(X, \varepsilon(X))$ is said
to be \emph{regular at infinity} (respectively, \emph{locally
compact at infinity}) if for every $E\in \varepsilon(X)$, there
exists $E' \in \varepsilon(X)$ such that $\overline{E'} \subset E$
($\overline{E'}$ is compact and $\overline{E'} \subset E$).
\end{dfn}

Obviously, locally compact at infinity implies regular at
infinity.

\begin{example}
As an example of a regular at infinity exterior space we can take
any Hausdorff locally compact space $X$ provided with its
cocompact externology. Now, if $K$ is a compact subset of $X$ and
we take the externology of open neighborhoods of $K$ in $X$, we
obtain a new exterior space which is  locally compact at infinity.
Next, we describe an exterior space which is not regular at
infinity: Consider the following planar differential system
$$\frac{d\alpha }{dt}=f(\alpha ,\beta ),
\hspace{10pt}\frac{d\beta}{dt}=af(\alpha ,\beta )$$ \noindent
where $f(\alpha ,\beta )>0,$ $f(\alpha +1,\beta )=f(\alpha
+1,\beta +1)=f(\alpha ,\beta+1)$ and $a$ is an irrational fixed
number. This system induces a flow on the torus $X=S^1\times S^1,$
which satisfies that each trajectory is dense in $X.$ Take the
externology $\varepsilon (X)$ constituted by those open subsets
$E$ such that for any $p\in X$ there is $r_0\in \mathbb{R}$ such
that $r\cdot p\in E,$ for all $r\geq r_0.$ Then it is easy to
check that the exterior space $(X,\varepsilon (X))$ is not regular
at infinity.
\end{example}

\begin{prop}\label{eleigualbarele}
If an exterior space $X=(X,\varepsilon(X))$ is regular at
infinity, then
\begin{itemize}
\item[(i)] $L(X)= \bar L(X)$
\item[(ii)]  $\F(X)=\bF(X)$, $e_0L(X)= \bar e_0\bar L(X)$
\item[(iii)]  $\C(X)=\bC(X)$, $e L(X)= \bar e \bar L(X)$.
\end{itemize}
\end{prop}

\begin{prop} \label{igual} Suppose that $X=(X,\varepsilon(X))$
is an exterior space.
\begin{itemize}
\item[(i)] If $X$ is locally path-connected, then  $\F(X)=\C(X)$,  $e_0L(X)= e L(X),$
\item[(ii)]  If $X$ is locally path-connected and regular at infinity, then
$\F(X)=\bF(X)=\C(X)=\bC(X)$, $e_0L(X)= \bar e_0\bar L(X)=e L(X)=
\bar e \bar L(X)$.
\end{itemize}
\end{prop}

\begin{prop}\label{masigual}
If an exterior space $X=(X,\varepsilon(X))$ is  locally compact at
infinity, then $L(X)=\bar L(X)$ is compact.
\end{prop}
\medskip
\begin{proof} Since $X$ is regular at infinity, by Proposition \ref{eleigualbarele}, $L(X) = \bar L (X)$.
Take $E_0\in \varepsilon(X)$ such that $\overline{E_0}$ is
compact, then the closed subset $\bar L (X)$ satisfies $\bar L (X)\subset
\overline{E_0}$. Therefore  $L(X)=\bar L(X)$ is compact.
\end{proof}

\begin{thm}\label{masmasigual} Let $X=(X,\varepsilon(X))$ be an exterior space and
suppose that $X$ is locally path-connected and locally compact at
infinity. Then,
\begin{itemize}
\item[(i)]   $L(X)=\bar L(X)$ is compact,
\item[(ii)] $e_0L(X)= \bar e_0\bar L(X)=e L(X)= \bar e
\bar L(X)=\F(X)=\bF(X)=\C(X)=\bC(X)$(any  end point is
representable by a point of the limit),
\item [(iii)] $\F(X)=\bF(X)=\C(X)=\bC(X)$ is a profinite compact space,
\item [(iv)] If $a \in \F(X)=\bF(X)=\C(X)=\bC(X)$, then $L_a(X)=\bar
L_a(X)=L_a^0(X)=\bar L_a^0(X)$ is a non-empty continuum.
\end{itemize}
\end{thm}

\medskip
\begin{proof}
As a consequence of  Propositions \ref{igual} and \ref{masigual},
it follows (i),  $\F(X)=\bF(X)=\C(X)=\bC(X)$ and $e_0L(X)= \bar
e_0 \bar L(X)=e L(X)= \bar e \bar L(X)$.

Now take $E_0\in \varepsilon(X)$ such that $\overline{E_0}$ is
compact, then $\bar \varepsilon'(X)= \{\overline{E}| E\in
\varepsilon(X), \overline{E}\subset \overline{E_0}\}$ is cofinal
and we have $\bar L(X) = \cap_{ \overline{E} \in \bar
\varepsilon'(X)} \overline{E}$ and $\bC(X)=
\lim_{\overline{E}\in\bar \varepsilon'(X)} c(\overline{E})$.

Note that any  end $a$ can be represented by $\{F\}_{F\in
c(\overline{E}),\overline{E}\in\bar \varepsilon'(X)}.$ Since $F$
is a non-empty component of $\overline{E}\subset \overline{E_0}$,
it follows that $F$ is closed ($F$ is a continuum). We also have
that the family of closed subset  $\{F\}_{F\in
c(\overline{E}),\overline{E}\in\bar \varepsilon'(X)}$ satisfies
the finite intersection property. Since $\overline{E_0}$ is
compact, one has that $L_a(X)=\bigcap_{F\in
c(\overline{E}),\overline{E}\in\bar \varepsilon'(X)}F$ is a
non-empty continuum (see Theorem 6.1.20 in \cite{Engelking}).
Therefore this end is representable by points of the limit space.
Since the map $e\colon L(X) \to \F(X)$ is continuous and $L(X)$ is
compact it follows that $\F(X)$ is compact. Moreover, since
$\F(X)\cong\lim \pi_0(E)$ is prodiscrete, taking into account that
$\F(X)$ is compact, we have that $\F(X)$ is a profinite compact
space.
\end{proof}

\begin{remark}
For an ANR exterior space $X,$ under some topological conditions
the shape of the limit space is determined by the resolution
$\varepsilon (X).$ Some applications of shape theory to dynamical
systems can be seen in \cite{Sanjurjo, Gabites}.
\end{remark}

\section{The category of $\r$-exterior flows}

We are going to consider the exterior space $\R^{\r}=(\R, \r)$,
where $\r$ is the following externology:
$${\r}=\{ U| U\hspace{3pt}\mbox{is open and there is}\hspace{3pt}
n\in \N\hspace{3pt}\mbox{such that} \hspace{3pt}(n, +\infty)
\subset U\}.$$ Note that a base for ${\r}$ is given by
$\mathcal{B}({\r})=\{(n, +\infty) | n \in \N\}.$

The exterior space $\R^{\r}$ plays an important role in the
definition of $\r$-\emph{exterior flow} below. Such notion mixes
the structures of dynamical system and exterior space:

\begin{dfn}  Let $M$ be an exterior space, $M_{\t}$ denote
the subjacent topological space and $M_{\d}$ denote the set $M$
provided with the discrete topology. An  $\r$-\emph{exterior flow}
is a continuous flow $\varphi \colon \R {\times} M_{\t} \to
M_{\t}$ such that $\varphi \colon \R^{\r}\bar\times M_{\d} \to  M$
is exterior and for any $t\in \R$, $F_t\colon M \bar \times I \to
M$, $F_t(x,s)=\varphi (ts,x)$, $s\in I,x \in M$, is also exterior.

An  $\r$-\emph{exterior flow morphism} of $\r$-exterior flows
$f\colon M \to N$ is a flow morphism  such that $f$ is exterior.
We will denote by $\bf{E^{\r} F}$ the category of $\r$-exterior
flows and $\r$-exterior flow morphisms.
\end{dfn}

Given an $\r$-exterior flow $(M,\varphi) \in  \bf{E^{\r} F}$, one
also has a flow $(M_{\t},\varphi) \in \bf{F}.$ This gives a
forgetful functor $$(\cdot)_{\t} \colon \bf{E^{\r} F} \to
\bf{F}.$$

Now given a continuous flow  $X=(X,\varphi)$, an open $N \in \t_X$
is said to be ${\r}$-\emph{exterior} if for any $x \in X$ there is
$T^x \in \r$ such that $\varphi (T^x \times \{x\}) \subset N.$  It
is easy to check that the family of $\r$-exterior subsets of $X$
is an externology, denoted by $\varepsilon^{\r}(X),$ which gives
an exterior space $X^{\r}=(X, \varepsilon^{\r}(X))$  such that
$\varphi \colon \R^{\r}\bar \times X_{\d} \to X^{\r} $ is exterior
and $F_t \colon X^{\r} \bar \times I \to  X^{\r}$,
$F_{t}(x,s)=\varphi( ts, x)$,  is also exterior for every
$t\in\R$. Therefore $(X^{\r},\varphi)$ is an $\r$-exterior flow
which   is said to be the $\r$-\emph{exterior flow associated to}
$X.$ When there is no possibility of confusion, $(X^{\r}, \varphi
)$ will be briefly denoted by $X^{\r}$.  Then we have a functor
$$(\cdot)^{{\r}} \colon \bf { F}\to \bf{E^{\r} F}.$$

Note that for a flow $(X, \varphi)$, if $E$ is an open subset such
that $\overline{E}$ is compact, then $E$  is an $\r$-exterior
subset  if and only if $\overline{E}$ is an ``absorbing region" in
the sense of Definition 1.4.2 in \cite{AH2003}. On the other hand,
the forgetful functor and the given constructions of exterior
flows are related as follows:

\begin{prop} The functor $(\cdot)^{{\r}} \colon \bf{ F}\to \bf{E^{\r} F}$
is left adjoint to the forgetful functor $(\cdot)_{\t} \colon \bf{E^{\r}
F}\to \bf{F}.$ Moreover $(\cdot)_{\t} \, (\cdot)^{{\r}}={\rm id}$
and $\bf{F}$ can be considered as a full subcategory of
$\bf{E^{\r} F}$ via $(\cdot)^{\r}.$
\end{prop}
\medskip
\begin{proof} Let $X$ be in $\bf { F}$ and $M$ be in $\bf{E^{\r} F}.$
If $f\colon X^{\r}\to M$ is a morphism in $\bf{E^{\r} F},$ then it
is clear that $f\colon X=(X^{\r})_{\t} \to M_{\t}$ is a morphism
in ${\bf F}.$ Now take $g\colon X \to M_{\t}$ a morphism in ${\bf
F}$ and $E\in \varepsilon (M).$ Given any $x\in X$ one has
$g(x)\in M$ and, taking into account that $M$ is an $\r$-exterior
flow, there exists $T^{g(x)}$ such that $T^{g(x)}\cdot g(x)\subset
E.$ This implies that $T^{g(x)}\cdot x \subset g^{-1}(E).$
Therefore $g^{-1}(E)\in \varepsilon(X^{\r})=\varepsilon^{\r}(X).$
\end{proof}

\subsection{End Spaces and Limit Spaces of an exterior flow}

In section \ref{tres} we have defined the end and limit spaces of
an exterior space. In particular, since any $\r$-exterior flow $X$
is an exterior space, we can consider the end spaces $\F(X),$
$\C(X),$ $\bF(X),$ $\bC(X)$ and the limit spaces $L(X),$ $\bar L
(X).$ Notice that one has the following properties:

\begin{prop} \label{invariantes} Suppose that $X=(X,\varphi)$ is an
$\r$-exterior flow. Then
\begin{itemize}
\item[(i)] The spaces  $L(X), \bar L  (X)$ are invariant;
\item[(ii)] There are trivial  flows induced  on $\F(X),$ $\C(X),$
$\bF(X)$ and $\bC(X).$
\end{itemize}
\end{prop}
\medskip
\begin{proof}(i): We have that $L(X)=\cap_{E\in \varepsilon(X)}E.$ Note that
for any $s\in \R,$ $\varphi_s(E)\in \varepsilon(X)$ if and only if
$E\in \varepsilon(X).$ Then $\varphi_s(L(X))=\varphi_s( \cap_{E\in
\varepsilon(X)}E)= \cap_{E\in
\varepsilon(X)}\varphi_s(E)=\cap_{E\in \varepsilon(X)}E=L(X).$ In
a similar way, it can be checked that $\bar L(X)$ is also
invariant.

(ii): For any $s\in \R$, consider the exterior homotopy
$F_{s}\colon X\bar \times I \to X$, $F_{s}(x,t)= \varphi(t s,x)$,
from $\id_X$ to $\varphi_s$. By Lemma \ref{lemahomotopia}, it
follows that $\id = \F(\varphi_s)$. Therefore the induced action
is trivial. In the other cases the proof is similar using Lemma
\ref{lemabarhomotopia}.
\end{proof}

As a consequence of this result, one has functors
 $L, \F, \C, \bar L, \bF, \bC\colon {\bf E^{\r}F} \to {\bf F}.$

\begin{prop}\label{homotopy} The functors
$L, \F, \C, \bar L, \bF, \bC\colon {\bf E^{\r}F} \to {\bf F}$
induce functors
$$L, \bar L\colon \pi {\bf E^{\r}F} \to \pi{\bf F},\hspace{10pt}\F,\C,\bF, \bC \colon
\pi{\bf E^{\r}F} \to {\bf F},$$ where the homotopy categories are
constructed in a canonical way.
\end{prop}

\subsection{The end point of a trajectory and the induced decompositions of
an exterior flow}

For a $\r$ exterior flow $X$, one has that  each trajectory has an
end point given as follows: Given $p \in X$ and $E \in
\varepsilon(X),$ there is $T^p\in \r$ such that $T^p \cdot p
\subset E.$ We can suppose that $T^p$ is path-connected and
therefore so is $T^p \cdot p;$ this way there is a unique
path-component $\omega_{\r}^0(p,E)$ (component $\omega_{\r}(p,
E)$) of $E$ such that $T^p \cdot p \subset \omega_{\r}^0(p, E)
\subset E$ ($T^p \cdot p \subset \omega_{\r}(p, E) \subset E$).
This gives set maps $\omega_{\r}^0 (\cdot , E) \colon X\to
\pi_0(E)$ and $\omega_{\r}^0 \colon X \to \F(X)$ ($\omega_{\r}
(\cdot , E) \colon X\to c(E)$ and $\omega_{\r} \colon X \to
\C(X)$) such that the following diagram commutes:
$$\xymatrix{& L(X)\ar[ldd]_{e_0} \ar[rdd]^{e}  \ar[d] & \\
&X  \ar[ld]^{\omega_{\r}^0} \ar[dr]_{\omega_{\r} }  & \\
 {\F(X)}\ar[rr]^{\theta} & & \C(X)\\
}$$

These maps permit us to divide a flow in simpler flows.

\begin{dfn} Let $X$ be an  $\r$-exterior flow. We will
consider $X_{(\r,a)}^0=(\omega_{\r}^0)^{-1}(a),$ $a\in \F(X)$ and
$X_{(\r,b)}=\omega_{\r}^{-1}(b),$ $b\in \C(X).$ The invariant
spaces $X_{(\r,a)}^0$  and $X_{(\r,b)}$ will be called the
$\omega_{\r}^0$-\emph{basin at} $a\in \F(X)$  and the
$\omega_{\r}$-\emph{basin at} $b\in \C(X),$ respectively.

The maps $\omega_{\r}^0$ and $\omega_{\r}$ induce the following
partitions of $X$ in simpler flows
$$X=\bigsqcup_{a\in \F(X) }X_{(\r, a)}^0,
\quad X=\bigsqcup_{b\in \C(X) }X_{(\r, b)} $$ that will be called
respectively, the $\omega_{\r}^0$-\emph{decomposition} and the
$\omega_{\r}$-\emph{decomposition} of the $\r$-exterior flow $X.$
\end{dfn}

Similarly, given $p\in X,$ if $\omega_{\r}^0(p,E)$ is the
path-component of $E$ such that $T^p \cdot p \subset
\omega_{\r}^0(p,E) \subset E,$ then we also have that $T^p \cdot p
\subset  \omega_{\r}^0(p, E) \subset \bar \omega_{\r}^0(p,
\overline{E})\subset \overline{E}$, where $ \bar \omega_{\r}^0(p,
\overline{E})$ is the unique path-component of $\overline{E}$
containing $T^p \cdot p.$ In the same way as above, we have maps
$\bar \omega_{\r}^0(\cdot , \overline{E}) \colon X\to
\pi_0(\overline{E})$ and  $\bar \omega_{\r}^0 \colon X \to
\bF(X).$ Analogously, we obtain set maps $\bar \omega_{\r}(\cdot ,
\overline{E}) \colon X\to c(\overline{E})$ and  $\bar \omega_{\r}
\colon X \to \bC(X)$ such that the following diagram commutes:

$$\xymatrix{& \bar L(X)\ar[ldd]_{\bar e_0} \ar[rdd]^{\bar e}  \ar[d] & \\
&X  \ar[ld]^{\bar \omega_{\r}^0} \ar[dr]_{\bar \omega_{\r} }  & \\
 {\bF(X)}\ar[rr]^{\bar \theta} & & \bC(X)\\
}$$

\begin{remark}
It is important to note that the maps $\omega_{\r}^0$,
$\omega_{\r}$ , $\bar \omega_{\r}^0$,  $\bar \omega_{\r}$ need not
be continuous.\end{remark}

As in the case above we can consider the corresponding $\bar
\omega_{\r}^0$-basin and  $\bar \omega_{\r}$-basin, denoted by
$\bar X_{(\r,a)}^0=(\bar \omega_{\r}^0)^{-1}(a)$, $\bar
X_{(\r,a)}= \bar \omega_{\r}^{-1}(a),$ respectively, and their
induced decompositions. We also note that the following diagram
commutes:
$$
\xymatrix{&X\ar[ld]_{ \omega_{\r}^0 }\ar[rd]^{\bar \omega_{\r}^0 }
\ar[ldd]^(.6){\omega_{\r} }\ar[rdd]_(.6){\bar \omega_{\r}}& \\
\F(X)\ar[d]\ar[rr]&&\bF(X)\ar[d]\\
\C(X)\ar[rr]&&\bC(X) }
$$

\begin{dfn} Let $X$ be an  $\r$-exterior flow. An end point $a\in \F(X)$
is said to be ${ \omega_{\r}^0 }$-\emph{representable} (similarly
for ${\omega_{\r} }, {\bar \omega_{\r}^0 },$ ${\bar \omega_{\r}}$)
if there is $p\in X$ such that $ \omega_{\r}^0 (p)=a.$
\end{dfn}

Denote by $ \omega_{\r}^0 (X)$ the space of
$\omega_{\r}^0$-representable end points (similarly,
$\omega_{\r}(X),$ $\bar \omega_{\r}^0 (X),$ $\bar \omega_{\r}
(X)$). Since the $\omega$-decompositions of $X$ are compatible
with the $e$-de\-com\-positions of the limit subspace, we have the
following commutative diagram of representable end points:
$$
\xymatrix{e_0 L(X) \ar[rr] \ar[rd] \ar[dd]& & eL(X) \ar[rd] \ar'[d][dd] & \\
&\omega_{\r}^0 (X)\ar[rr] \ar[dd] & &\omega_{\r} (X) \ar[dd]\\
\bar e_0 \bar L(X) \ar'[r][rr] \ar[rd]& & \bar e\bar L(X) \ar[rd]& \\
&\bar \omega_{\r}^0 (X)\ar[rr] & &\bar \omega_{\r} (X)\\
}
$$

\begin{prop}\label{erepresentableuno}
Let $X$ be an $\r$-exterior flow. Then
\begin{itemize}
\item[(i)]  $\omega^{\r}(p) \subset \bar L_{\bar
\omega_{\r}(p)}(X),$ for any $p\in X.$
\item[(ii)]
If $a\in \bC(X)$ is $\bar\omega_{\r}$-representable (that is,
$\bar X_{(\r, a)}=\omega_{\r}^{-1}(a)\not = \emptyset$) and there
exists $p\in\bar X_{(\r, a)}$ such that $\omega^{\r}(p) \not =
\emptyset$, then $a$ is $\bar e$-representable.
 \end{itemize}
\end{prop}
\medskip
\begin{proof} (i) If $q\in \omega^{\r}(p),$ then
$q\in \cap_{T\in \r}\overline{T \cdot p}$. On the other hand,
given $E\in \varepsilon(X),$ if $\bar \omega_{\r}(a,\overline{E})$
is the connected component of $\overline{E}$ determined by $a$,
since $\bar\omega_{\r}(p)=a$, there is $T\in \r$ such that
$\overline{T \cdot p}\subset \bar \omega_{\r}(a,\overline{E}).$
Then $q\in \bar\omega_{\r}(a,\overline{E})$ for every $E\in
\varepsilon(X)$ and $q\in \bar L(X).$ This implies that $\bar
e(q)=a$ for any $q\in \omega^{\r}(p).$ Therefore $\omega^{\r}(p)
\subset \bar L_{\bar \omega_{\r}(p)}(X)$. (ii) follows from
(i).
\end{proof}

\begin{prop}\label{erepresentabledos} Let $X$ be an $\r$-exterior flow
and for  $p\in X$ denote $\gamma_{\r}(p)=\{r\cdot p|r \ge 0\}.$ If
$\overline{\gamma_{\r}(p)}\cap \bar L(X) \not = \emptyset,$ then
$\bar \omega_{\r}(p)$ is $\bar e$-representable.
\end{prop}
\medskip
\begin{proof} We note that  $\overline{\gamma_{\r}(p)}=\gamma_{\r}(p)
\cup \omega^{\r}(p).$ If $p\in \bar L (X)$, then $\bar
\omega_{\r}(p)= \bar e (p)$ and $\bar \omega_{\r}(p)$ is $\bar
e$-representable. If $p \not \in \bar L(X),$ then $\gamma_{\r}(p)
\cap \bar L (X)=\emptyset$ and $\overline{\gamma_{\r}(p)}\cap \bar
L(X) =  \omega^{\r}(p).$ Then  $\omega^{\r}(p) \not = \emptyset$
and, by Proposition \ref{erepresentableuno} above, $\bar
\omega_{\r}(p)$ is $\bar e$-representable.
\end{proof}

\section{End and Limit spaces of a flow via exterior flows}

Recall that we have considered the functor:
$$(\cdot)^{{\r}} \colon \bf { F}\to \bf{E^{\r} F}$$
and the functors: $$L, \F, \C, \bar L, \bF, \bC\colon {\bf
E^{\r}F} \to {\bf F}.$$ Therefore we can consider the composites:
$$L^{{\r}}=L(\cdot)^{{\r}},\hspace{5pt}\F^{{\r}}=\F(\cdot)^{{\r}},\hspace{3pt}
\C^{{\r}}=\C(\cdot)^{{\r}},\hspace{5pt}\bar L^{{\r}}=\bar
L(\cdot)^{{\r}},\hspace{5pt}\bF^{\r}=\bF(\cdot)^{{\r}},\hspace{5pt}
\bC^{{\r}}=\bC(\cdot)^{{\r}}$$ to obtain functors $L^{{\r}},$
$\F^{{\r}},$ $\C^{{\r}},$ $\bar L^{{\r}},$ $\bF^{\r},$
$\bC^{{\r}}\colon {\bf F} \to {\bf F}.$

In this way, given a flow $X$, we have the end spaces
$\F^{\r}(X)=\F(X^{\r})$,  $\C^{\r}(X)=\C(X^{\r})$, the limit space
$L^{\r}(X)=L(X^{\r})$, the bar-end spaces
$\bF^{\r}(X)=\bF(X^{\r})$, $\bC^{\r}(X)=\bC(X^{\r})$ and the
bar-limit space $\bar L^{\r}(X)=\bar L(X^{\r}).$

Similarly, using the associated exterior flow $X^{\r},$ we denote
$$X_{(\r,a)}^0=(\omega_{\r}^0)^{-1}(a), \quad a \in  \F^{\r}(X)$$
$$X_{(\r,a)}=\omega_{\r}^{-1}(a), \quad a \in  \C^{\r}(X)$$
The maps  $\omega_{\r}^0$,  $\omega_{\r}$  induce the following
partitions of $X$ in simpler flows
$$X=\bigsqcup_{a\in \F^{\r}(X) }X_{(\r, a)}^0,
\quad X=\bigsqcup_{a\in \C^{\r}(X) }X_{(\r, a)} $$ that  will be
called respectively, the $\omega_{\r}^0$-decomposition and the
$\omega_{\r}$-decomposition of the flow $X.$ Now take
 $$\bar X_{(\r,a)}^0=(\bar \omega_{\r}^0)^{-1}(a), \quad a \in   \bF^{\r}(X)$$
$$\bar X_{(\r,a)}=\bar \omega_{\r}^{-1}(a), \quad a \in   \bC^{\r}(X)$$
the maps  $\bar \omega_{\r}^0$,  $\bar \omega_{\r}$  induce the
$\bar \omega_{\r}^0$-decomposition and the  $\bar
\omega_{\r}$-decomposition of the flow $X:$
$$X=\bigsqcup_{a\in \bF^{\r}(X) }\bar X_{(\r, a)}^0,
\quad X=\bigsqcup_{a\in \bC^{\r}(X) }\bar X_{(\r, a)}.$$

\bigskip It is interesting to consider the following equivalence of
categories: Given any flow $\varphi \colon \R \times X \to X,$ one
can consider the \emph{reversed flow} $\varphi' \colon \R \times X
\to X$ defined by $\varphi'( r, x)=\varphi (-r,x),$ for every
$(r,x) \in \R \times X.$ The correspondence, $(X,\varphi)\to
(X,\varphi'),$ gives rise to a functor
$$(\cdot)' \colon {\bf F} \to {\bf F}$$ which is an equivalence of
categories and verifies $(\cdot)' (\cdot)'= \id.$ Using the
composites
$$L^{{\l}}=(\cdot)'L^{{\r}}(\cdot)',\hspace{5pt}
\F^{{\l}}=(\cdot)'\F^{{\r}}(\cdot)',\hspace{5pt}\C^{{\l}}=
(\cdot)'\C^{{\r}}(\cdot)'$$ $$\hspace{5pt}\bar
L^{{\l}}=(\cdot)'\bar L^{{\r}}(\cdot)',
\hspace{5pt}\bF^{\l}=(\cdot)'\bF^{\r}(\cdot)',\hspace{5pt}
\bC^{{\l}}=(\cdot)'\bC^{{\r}}(\cdot)'$$ we obtain the new functors
$L^{{\l}}, \F^{{\l}}, \C^{{\l}},\bar L^{{\l}}, \bF^{\l},
\bC^{{\l}}\colon {\bf F} \to {\bf F}$ and the decompositions
$$X=\bigsqcup_{a\in \F^{\l}(X) }X_{(\l, a)}^0, \quad X=
\bigsqcup_{a\in \C^{\l}(X) }X_{(\l, a)} $$
$$X=\bigsqcup_{a\in \bF^{\l}(X) }\bar X_{(\l, a)}^0,
\quad X=\bigsqcup_{a\in \bC^{\l}(X) }\bar X_{(\l, a)}.$$

\begin{remark} All decompositions above
can be considered as generalizations for a continuous flow of  the
disjoint union of ``stable" (or ``unstable" for the dual case)
submanifolds of a differentiable flow(see \cite{Salamon},
\cite{Sanjurjo}, \cite{Smale}).
\end{remark}

We  note that the decompositions of a flow $X$ are compatible with
decompositions of  limit subspaces.

 For a Morse function \cite{Milnor} $f\colon M \to \R$, where
$M$ is a compact $T_2$ Riemannian manifold, one has that the
opposite of the gradient of the $f$ induces a flow with a finite
number of critical points. In this case, we have that $M$ is
locally path-connected and the flow is $\r$-regular at infinity.
Then we have all the properties obtained by Theorem
\ref{masmasigual}. For instance we can take the height function of
a 2-torus:

\begin{example}
Let $\varphi \colon \R \times (S^1\times S^1) \to S^1\times S^1$
be the flow induced by the opposite of the gradient of the height
function with four critical points:
\begin{center}
\includegraphics[scale=0.3]{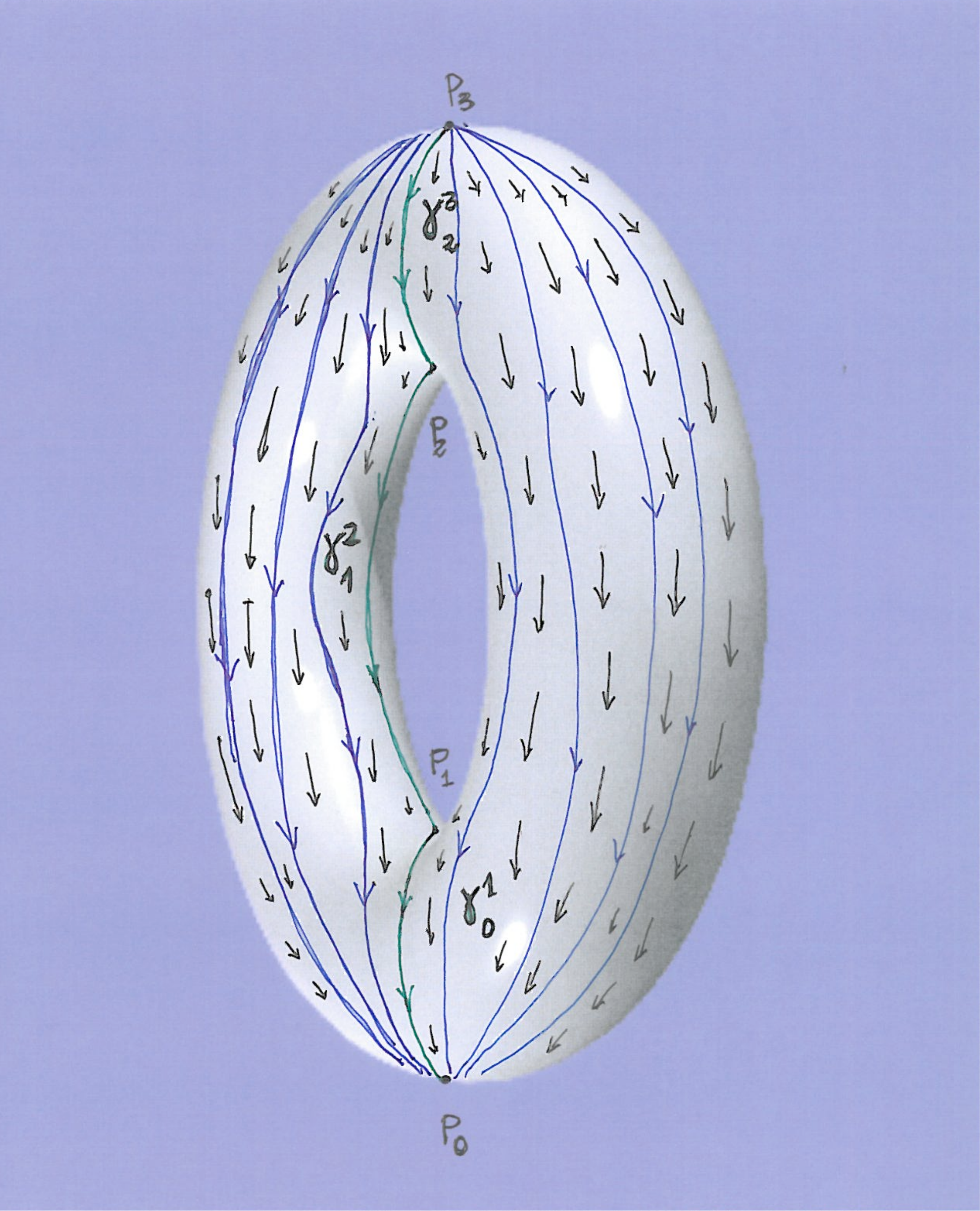}
\end{center}
In this example, the limit space, end space and decomposition of
the flow and the reverse flow are given in the following table:
\begin{flushleft}
\begin{tabular}
[c]{| l | l | }\hline $L^{\r}=\{P_0,P_1,P_2,P_3\}$ &
$\F^{\r}=\{P_0,P_1,P_2,P_3\}$  \\\hline $L_{P_i}^{\r}=\{P_i\}$  &
\\\hline $X_{(\r,P_3)}=\{P_3\}$ &$X_{(\r,P_2)}=\{P_2\}\cup
\gamma_2^3 \cup{\tilde \gamma}_2^3 $  \\\hline
 $X_{(\r,P_1)}=\{P_1\}\cup \gamma_1^2 \cup{\tilde \gamma}_1^2 $ &
 $X_{(\r,P_0)}=(S^1\times S^1)\setminus \bigcup_{i=1}^3 X_{(\r,P_i)}$  \\\hline
$L^{\l}=\{P_0,P_1,P_2,P_3\}$\, ~~~~~&
$\F^{\l}=\{P_0,P_1,P_2,P_3\}~$~ \, ~~~\\\hline
$L_{P_i}^{\l}=\{P_i\}$  &  \\\hline $X_{(\l,P_0)}=\{P_0\}$
&$X_{(\l,P_1)}=\{P_1\}\cup \gamma_0^1 \cup{\tilde \gamma}_0^1 $
\\\hline
 $X_{(\l,P_2)}=\{P_2\}\cup \gamma_1^2 \cup{\tilde \gamma}_1^2 $ &
 $X_{(\l,P_3)}=(S^1\times S^1)\setminus ~\bigcup_{i=0}^2 X_{(\l,P_i)}$\\\hline
\end{tabular}
\end{flushleft}
\end{example}

Now we consider a flow induced by a lineal differential equation
on $\R^2$ that also induces a new flow on the Alexandrov one-point
compactification $S^2=\R^2 \cup \{\infty\}.$

\begin{example} Consider on  $S^2$ the flow induced by
$ \varphi(r, (u_1,u_2))= (e^{r\lambda_1}u_1, e^{r\lambda_2}u_2),$
$u_1,u_2 \in  \R$, $\varphi(r, \infty)= \infty$, ($\lambda_1>0,$
$\lambda_2<0$)

\begin{center}
\includegraphics[scale=0.7]{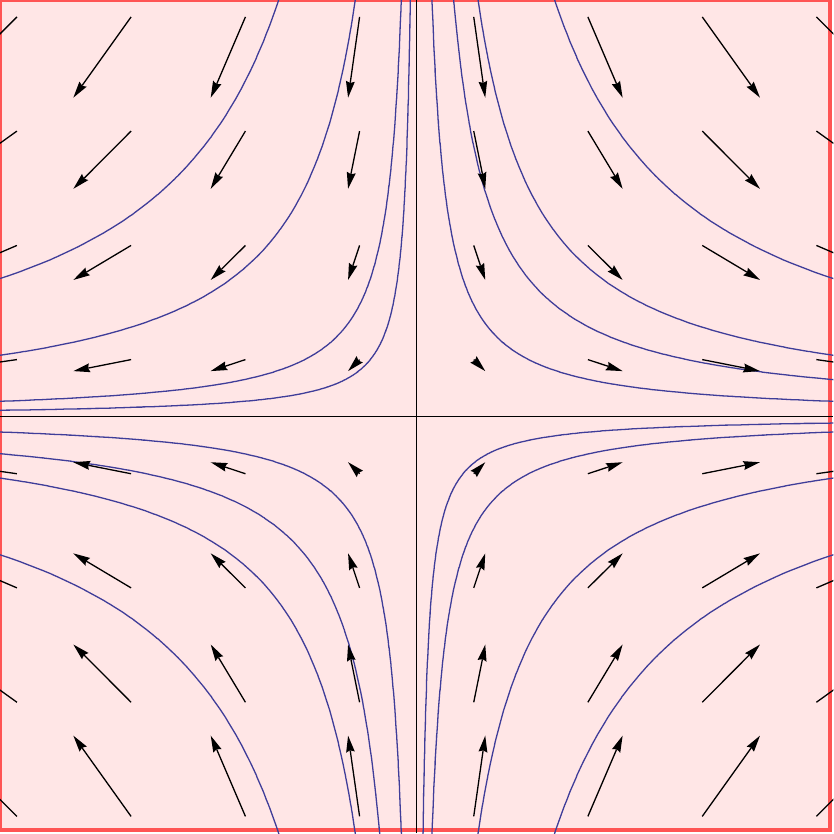}
\end{center}

The limit spaces, end spaces and decomposition of the flow $(S^2,
\phi)$ (as well as the reverse flow) are given in the following
table:

\bigskip
{
\begin{tabular}
[c]{| l | l | }\hline
  $L^{\r}=\{0, \infty\}$  & $\F^{\r}=\{0, \infty\}$  \\\hline
$L_0^{\r}=\{0\}$  &$L_{\infty}^{\r}=\{ \infty\}$  \\\hline
 $X_{(\r,0)}=\{0\} \times \R$  &$X_{(\r, \infty)}=((\R\setminus\{0\})\times \R)\cup \{\infty\}$ \\\hline
  $L^{\l}=\{0, \infty\}$  & $\F^{\l}=\{0, \infty\}$  \\\hline
  $L_0^{\l}=\{0\}$  &$L_{\infty}^{\l}=\{ \infty\}$ \\\hline
$X_{(\l,0)}=\R\times \{0\}$  & $X_{(\l, \infty)}=(\R \times
(\R\setminus\{0\}))\cup  \{\infty\}$    \\\hline

\end{tabular}
}
\end{example}

\section{Relations between limit and end spaces of a flow  and its dynamical properties}

\subsection{Periodic points}

The relation of the limit space of a flow or an $\r$-exterior flow
and the subflow of periodic points is analysed in the following
results:

\begin{lem}\label{uno} If $X$ is an $\r$-exterior flow, then
$P(X)\subset L(X)$. In particular, if $X$ is a flow, then
$P(X)\subset L^{\r}(X).$
\end{lem}
\medskip
\begin{proof} Take $x$ a periodic point and $E\in \varepsilon(X)$ arbitrary.
Then there exists $T\in \r$ such that $T\cdot x \subset E.$ Since
$x$ is periodic, $T\cdot x= \R \cdot x$ and taking into account
that $x\in \R \cdot x,$ we have that that $x\in E.$
\end{proof}

\begin{lem}\label{dos}
Let $X$ be a flow and suppose that $X$ is a $T_1$-space. Then, for
every $x\in X$  the following statements are equivalent:
 \begin{itemize}
  \item[(i)] $x$ is a non-periodic point.
  \item[(ii)]  $X\setminus \{x\}$ is  an $\r$-exterior subset of $X.$
\end{itemize}
\end{lem}
\medskip
\begin{proof} (i) implies (ii): Take $y\in X$; if the trajectory of
$y$ is different of the trajectory of $x,$ then for every $T\in
\r,$ $T\cdot y\subset X\setminus \{x\}.$ If $y$ is in the
trajectory of $x,$ considering that $x$ is not periodic, one can
find $T\in \r $ such that $T\cdot y \subset X\setminus \{x\}.$
Then, one has that $X\setminus \{x\} \in \varepsilon^{\r} (X).$
Conversely, suppose that $x$ is a periodic point, by Lemma
\ref{uno} above $X\setminus \{x\}$ is not $\r$-exterior.
\end{proof}

\begin{thm}\label{periodic} Let $X$ be a flow and suppose that $X$ is a
$T_1$-space. Then $$P(X)=L^{\r}(X).$$
\end{thm}
\medskip
\begin{proof} Let $x\in X\setminus P(X)$.
Then, by Lemma \ref{dos}, one has that $X\setminus \{x\} \in
\varepsilon^{\r} (X)$ and $$P(X)=X \setminus (\bigcup_{x \not \in
P(X)}\{x\})=\bigcap_{x \not \in P(X)} X\setminus \{x\}\supset
\bigcap_{E\in \varepsilon^{\r}(X)}E=L^{\r}(X).$$ Now the results
follows from Lemma \ref{uno}.
\end{proof}

Taking into account the result above, if $X$ is flow and $X$ is a
$T_1$ space we also have that
$$L^{\r}(X)=P(X)\subset P^{\r}(X)\subset \Omega^{\r}(X) \subset X.$$

\subsection{Limit spaces and invariant sets}

\begin{lem}Given a flow $\varphi\colon \R \times X \to X$
and $A \subset X$ we have that
$$\inv(A)= \bigcap_{r\in \R} \varphi_r (A).$$
\end{lem}
\medskip
\begin{proof}  If $p \in \inv(A),$ then $\R\cdot p \subset A.$
Notice that $p=\varphi_r(\varphi_{-r}(p))\in \varphi_r(A)$ so $p
\in  \bigcap_{r\in \R} \varphi_r (A).$ Conversely, if $p \in
\bigcap_{r\in \R} \varphi_r (A),$ then $p=\varphi_r(x_r),$ $x_r
\in A.$ This implies that $\varphi_{-r} (p)=x_r \in A$ and $p \in
\inv (A).$
\end{proof}

Then, by Proposition \ref{invariantes}, we obtain:

\begin{prop} If $X$ is an $\r$-exterior flow, then
\begin{itemize}
\item[(i)] $L(X)=\lim_{E\in \varepsilon(X)} E=\lim_{E\in \varepsilon(X)} \inv (E),$
\item[(ii)] $\bar L(X)=\lim_{E\in \varepsilon(X)} \overline{E}=
\lim_{E\in \varepsilon(X)} \inv (\overline{E}).$
\end{itemize}
\end{prop}

Next we give a characterization of the points lying in the
difference $\bar L(X) \setminus  L(X):$

\begin{prop}Let $X$ be an $\r$-exterior flow and $x\in \bar L(X)$.
Then $x \in\bar L(X) \setminus  L(X)$ if and only if there exist
$E\in \varepsilon(X)$ and  $t\in \R$ such that $t \cdot x \in \Fr
(E)$.
\end{prop}

Note that propositions above can be applied to $L^{\r}(X)$ and
${\bar L}^{\r}(X)$ for a continuous  flow $X$.

\subsection{Limits and $\Omega $-limits}

In the following result, we analyse the relationship between the
$\Omega^{\r}$-limit and the bar-limit induced by an externology.

\begin{lem}\label{omega_bar_limit} If $X$ is an $\r$-exterior flow, then
$$\Omega^{\r}(X)\subset \bar L (X)$$
\end{lem}
\medskip
\begin{proof} If $E \in \varepsilon (X)$, for every $x \in X$ there exist
$T\in \r$ such that $T\cdot x \subset E.$ Therefore
$\overline{T\cdot x} \subset \overline{E}.$ By definition this
implies that $\omega^{\r}(x) \subset \bar L (X)$ for every $x\in
X.$ Hence $\Omega^{\r}(X)\subset \bar L(X).$
\end{proof}

\begin{prop}\label{exterior_bar_limit}
Let $X$ be an $\r$-exterior flow, $\varepsilon(X)$ its externology
and $x\in X.$ Then there exists $V_x \in {({\t}_X)}_x$ such that
$X\setminus \overline{V}_x \in \varepsilon(X)$ if and only if $x
\not \in \bar L (X)$.
\end{prop}
\medskip
\begin{proof} Suppose that $X\setminus \overline{V}_x  \in \varepsilon(X)$.
Since $V_x \cap (X\setminus \overline{V}_x )= \emptyset $, then $x
\not \in \overline{X\setminus \overline{V}_x}$. Since $X\setminus
\overline{V}_x\in \varepsilon(X)$, it follows that $x\not \in
\bigcap_{E \in \varepsilon(X)} \overline{E}=\bar L (X)$.

Conversely, if $x \not \in \bar L (X),$ then there is $E \in
\varepsilon(X)$ such that $x\not \in \overline{E}.$ Taking $V_x=
X\setminus \overline{E}=$int$(X\setminus E)$, we have that
$X\setminus \overline{V}_x=$int$(X\setminus
V_x)=$int$(\overline{E}) \supset E.$ Consequently, $X\setminus
\overline{V}_x\in \varepsilon(X)$.
\end{proof}



\begin{cor}
Let $X$ be an $\r$-exterior flow, $\varepsilon(X)$ its externology
and $x\in X.$ If there exists $V_x \in {({\t}_X)}_x$ such that
$X\setminus \overline{V}_x\in \varepsilon(X)$, then $x \not \in
\overline{\Omega^{\r}(X)}.$
\end{cor}
\medskip
\begin{proof} It is a consequence of Proposition \ref{exterior_bar_limit}
and  Lemma \ref{omega_bar_limit}.
\end{proof}

\begin{cor}
Let $X$ be a flow  and $x\in X.$ If there exists $V_x \in
{({\t}_X)}_x$ such that $X\setminus \overline{V}_x$ is
$\r$-exterior, then $x\not \in \overline{\Omega^{\r}(X)}.$
\end{cor}

\begin{lem} Let $X$ be a flow and $X$ is a locally compact regular space.
If  $x \not \in\overline{ \Omega^{\r}(X)},$ then there exists $V_x
\in {({\t}_X)}_x$ such that  $X\setminus \overline{V}_x$ is
$\r$-exterior.
\end{lem}
\medskip
\begin{proof} Suppose that $x \not \in\overline{ \Omega^{\r}(X)}.$
Since $X$ is locally compact, there is a compact neighborhood $K$
at $x$ such that $K \cap \Omega^{\r}(X)= \emptyset.$ Take $y\in X$
and assume that for every $T \in \r,$ $T\cdot y \cap K \not =
\emptyset .$ Then there is a sequence $t_n\rightarrow +\infty $
such that $t_{n}\cdot y\in K.$ Since $K$ is compact, one can take
a subsequence $t_{n_i}\rightarrow +\infty$ such that $t_{n_i}\cdot
y \to u \in K.$ This fact implies that $u\in K \cap \omega^{\r}(y)
\subset K\cap \Omega^{\r}(X),$ which is a contradiction.
Therefore, there is $T$ such that $T\cdot y \cap K= \emptyset.$
Since $X$ is regular, one has that there is $V_x \in {({\t}_X)}_x$
such that $\overline{V}_x \subset K$ and $X\setminus
\overline{V}_x$ is  $\r$-exterior.
\end{proof}

\begin{cor} Let $X$ be a flow. If $X$ is a locally compact regular space,
then $\bar L^{\r}(X) \subset \overline{\Omega^{\r}(X)}.$
\end{cor}
\medskip
\begin{proof} If $x \not \in \overline{\Omega^{\r}(X)}$, by the lemma above,
there exists $V_x \in {({\t}_X)}_x$ such that  $X\setminus
\overline{V}_x$ is  $\r$-exterior. By Proposition
\ref{exterior_bar_limit}, it follows that $x \not \in \bar
L^{\r}(X).$
\end{proof}

\begin{thm} \label{main}Let $X$ be a flow. If $X$ is a locally compact
regular space, then $\bar L^{\r}(X) = \overline{\Omega^{\r}(X)}.$
\end{thm}
\medskip
\begin{proof} By the corollary above $ \bar L^{\r}(X)
\subset \overline{\Omega^{\r}(X)}$ and by Lemma
\ref{omega_bar_limit}, $\Omega^{\r}(X)\subset \bar L (X).$
\end{proof}

\begin{cor} \label{reticulo} Let $X$ be a flow. If $X$ is a locally compact
$T_3$  space, then $L^{\r}(X) = P(X)$, $ \bar L^{\r}(X)=
\overline{\Omega^{\r}(X)}$ and $$L^{\r}(X)=P(X)\subset
P^{\r}(X)\subset \Omega^{\r}(X) \subset
\overline{\Omega^{\r}(X)}=\bar L^{\r}(X).$$
\end{cor}
\medskip
\begin{proof}
This is a consequence of Theorems \ref{periodic} and \ref{main}.
\end{proof}

\end{document}